\documentstyle[12pt]{article}
\def\ni{\noindent}
\begin{document}
\begin{center}
{\large \bf Life and work of the Mathemagician Srinivasa Ramanujan}\\
\bigskip
{\bf K. Srinivasa Rao}\\                          
The Institute of Mathematical Sciences, Chennai – 600 113\\.
(E-mail : rao@imsc.ernet.in)\\
\end{center}
\bigskip

\ni {\bf Introduction}\\

Srinivasa Ramanujan, hailed as one of the greatest mathematicians of this century,
left behind an incredibly vast and formidable amount of original work, which has 
greatly influenced the development and growth of some of the best research work 
in mathematics of this century.  He was born at Erode, on Dec. 22, 1887. There 
were no portents to indicate that he would, in a short life-span of 32 years 4 
months and 4 days, become comparable to the all-time great Euler, Gauss and 
Jacobi, for {\it natural genius}.\\ 

There are two aspects of interest to biographers and mathematicians regarding
Ramanujan: his life and his work.  Mathematicians, who are interested in his 
work, have to contend with not only his publications in journals which are precise 
and profound, but also with his Notebooks which are a treasure house of 
intriguing results stated without proofs and lacking perspective with 
contemporary mathematical work. Those who attempt to write biographic articles 
on Ramanujan have to surmount the time barrier to reconstruct a story from all the 
indirect information accessible and to them, Hardy on Ramanujan [1] is akin to 
Boswell on Samuel Johnson. The challenge to the mathematicians who work on 
any of his thousands of recorded results, which are still shrouded in mystery, is to 
prove the same with what was accessible to Ramanujan in those days in the form 
of books and publications. While the individual writer's' perception of Ramanujan 
will depend upon his/her background and imagination, the task of  the 
mathematician is perhaps unenviable, in comparison.  \\

Anyone who ever heard of Srinivasa Ramanujan and reads the compelling 
rags-to-intellectual-riches story of Ramanujan contained in the two Notices, one by 
G.H. Hardy and the other by  Dewan Bahadur R. Ramachandra Rao and P.V. 
Seshu Iyer,  published in the {\it Collected papers of Srinivasa Ramanujan} [2],
would 
be moved by the achievements of the unorthodox mathematical genius under 
adverse circumstances. The lack of formal education, lack of appreciation and a 
job, in the beginning of his career and ill health during the last few years of his 
life, did not prevent him from being creative in Mathematics. This is indeed 
something not easy to comprehend, for often one would buckle under similar 
trying circumstances. In these lectures, I will present an account of his romantic 
life, provide a few glimpses into his mathematics and relate the increasing interest  
in his work and its relevance even today.\\

\medskip

\ni {\bf Formal education}\\

Ramanujan's' father, Mr. K. Srinivasa Iyengar, was an accountant to a cloth 
merchant in Kumbakonam. His mother was Komalattammal and Erode was her 
parental home. He was the first of three sons to his parents. Very little is known 
about his father and not even a  photograph  of his seems to be available. His 
mother was convinced of the greatness of Ramanujan  and she zealously protected 
and projected his interests all through his life. She is portrayed as a shrewd, 
cultured lady and her photograph is available in some books on  Ramanujan.\\

Ramanujan  was sent to Kangeyam Primary School in Kumbakonam at the age of 
seven. During his school days, he impressed his classmates, senior students and 
teachers with his extraordinary intuition and astounding proficiency in several 
branches of mathematics - viz. arithmetic, algebra, geometry, number theory and 
trigonometry. In later years a friend of his, C.V. Rajagopalachari, recounted the 
following incident ([3], p.83) 
which happened when Ramanujan was in his third form: In an arithmetic class on 
division, the teacher said that if three bananas were given to three boys, each boy 
would get a banana. The teacher generalised this idea and said that any number 
divided by itself would give one. Ramanujan asked: \\
        
{\it   Sir, if no banana is distributed to no student, will everyone still get a
banana ?}\\  

Another friend who took private tuition from Ramanujan also recalled [4] that
Ramanujan used to ask about the value of zero divided by zero and then answer 
that it can be anything since the zero of the denominator may be several times the 
zero of the numerator and vice versa and that the value cannot be determined. He 
stood first in the Tanjore District Primary Examinations held in November 1897, 
and this entitled him to a half-fee concession in the Town High School at 
Kumbakonam, where he studied from 1898 to 1903, until he passed the  
Matriculation Examination of the University of Madras (1904).\\

At the age of 12, Ramanujan is said to have worked out the properties of 
arithmetical, geometrical and harmonic progressions. Once a senior school 
student [3], posed to 
Ramanujan, who was in the fourth year at school, the 
following problem: 

$$If \ \surd{x}  + y = 7\quad {\rm and}\quad  x + \surd{y} = 11,\ what\ are\ the 
\ values\ of\ x\ and\ y\ ? $$

Ramanujan''s immediate reply to this questionn--– which was expected to be tackled
by only a sixth year student -- that $ x = 9$ and $y = 4$, won for him a friend who
in later years took him to the collector of Nellore\footnote{If one does not guess
this answer, the result can be obtained by setting $x = m^2, y = n^2$, then
take the difference between the two simultaneous equations and factorise to get:
$(m-n)(m+n-1) = 4$, which has integer solutions only for $m=3, n=2$ and hence 
$x=9, y=4$.}.\\

The senior mathematics teacher of the school, Ganapathy Subbier, had such
confidence in Ramanujan''s ability that year after year he entrusted Ramanujan with
the task of preparing a conflict free time-table [5] for the school, which had
about 1500 students and 30 or more teachers. Ramanujan won prizes for his
outstanding performance in mathematics and mastered Loney''s{\it  Trigonometry,
Part II},
in his fourth year at school. He won many prizes [6] in his second, fourth and
sixth years at High School.\\ 

To augment the family income, Ramanujan''s mother took in a couple of students from
Tirunelveli and Tiruchirapalli as boarders. Noticing Ramanujan's precocity in
mathematics these undergraduate students are purported to have given him an
elementary introduction to all branches of mathematics. In 1903, through these
friends from the Kumbakonam Government College, Ramanujan obtained G.S. Carr's':
{\it A Synopsis of Elementary Results, a book on Pure Mathematics}, 
which contained
propositions, formulae and methods of analysis with abridged demonstrations,
published in 1886.\\ 

Carr presented in this book 4865 formulae [7, p.3], without proofs, in
algebra, trigonometry,
analytical geometry and calculus. This book is similar to the modern day
compilations like the {\it Table of Integrals, Series, and Products}, by I.S.
Gradshteyn
and I.M. Ryzhik (Academic Press, New York, 1994). Prof. P.V. Seshu
Aiyar and Mr. R. Ramachandra Rao, in their biographies of Ramanujan [2] state that: 
\\ 

%\leftskip 1cm \rightskip 1cm

\ni {\it It was this book which awakened his genius. He set himself to establish 
    the formulae given therein. As he was without the aid of other books, each  
    solution was a piece of research so far as he was concerned.}\\

%\leftskip 0cm \rightskip 0cm

\ni It is the considered opinion of many (cf. Kanigel [8], p.57) that in proving
one 
formula, he discovered many others and thus, Ramanujan laid for himself a 
foundation for higher mathematics. Also,  at about this time, he started noting his 
results in Notebooks.\\

The first public recognition of his extraordinary prowess came when he was 
awarded a special prize --– the Sri K. Ranganatha Rao Prize  e--– at the annual prize 
distribution ceremony of the Town High School, in 1904, for proficiency in 
mathematics. Ramanujan passed his Matriculation Examination in 1904 and  
joined the Government Arts College in Kumbakonam. As a result of his success 
in a competitive examination in Mathematics and English composition, he 
secured the Junior Subrahmanyam Scholarship. In the F.A. (First Examination in 
Arts) Class, Ramanujan had to study English, Sanskrit, Mathematics, Physiology 
and the History of Rome and Greece. Partly due to his pre-occupation with 
researches into mathematics, he neglected the study of other subjects. He went  to 
his mathematics lecturer with a number of original and very ingenious results in 
finite and infinite series. Prof. P.V. Seshu Aiyar exhorted him but advised him not 
to neglect the study of other subjects. Unfortunately, he did not pass in English 
and Physiology and hence was not promoted to the senior F.A. class in January 
1905. He lost his scholarship. His mother, who played a domineering role in his 
life, tried to persuade the Principal of the Government Arts College to take note of  
Ramanujan's extraordinary mathematical ability and appealed for a continuance of 
the scholarship, but to no avail.\\

Ramanujan''s failure to get promoted to the senior F.A. class marked the 
beginning of a very trying period in his life. It is not clear what he did in 1905, 
when he discontinued  his studies and spent some months in (the present day) 
Andhra Pradesh region, when he set out from Kumbakonam, for the first time. He 
joined Pachaiyappa''s College in Madras, in the F.A. class again, in 1906. One of 
his classmates, T. Devaraja Mudaliar, ([9], p.63 and p.65)
recalls that the Chief Professor of 
Mathematics, P. Singaravelu 
Mudaliar, considered an acquisition by Pachaiyappa''s College since he had the
reputation of being a very successful teacher for the B.A. class, waited for 
Ramanujan''s assistance to solve difficult problems in mathematical journals. He  
also recalls that a junior mathematics teacher of the F.A. class, Prof. N. 
Ramanujachari, allowed Ramanujan to go to the board to show the solutions to 
the difficult problems in algebra or trigonometry using  fewer steps than the ones 
used by him. Senior students of the B.A. Class also sought Ramanujan''s help in
mathematics [10].\\

Ramanujan who was a strict vegetarian should have abhorred the dissection of the 
frog in the Physiology classes. Once, to a question on the digestive system, he is 
supposed to have provided a skimpy answer which he concluded with [11]: {\it" Sir, 
this is my undigested product of the Digestion chapter. Please excuse me'}'.
Another classmate of his at Pachaiyappa's College recalls [12] that Ramanujan
{\it "rarely got more than 10 % in Physiology, for which subject he had supreme 
contempt and got something more, say 15 \% to 20 \% in Greek and Roman 
History, but managed to get about 25 \% in English}." However, Ramanujan 
considered [12] {\it "the problems given in$\cdots$  textbooks in Geometry,
Algebra, and Trigonometry}" to
be {\it "mental  sums}.\\

In 1906, while studying at Pachaiyappa''s College, Ramanujan lived with his 
grandmother in a house in a lane in George Town, Madras. After about three 
months, Ramanujan  fell ill and discontinued his studies. However, he appeared 
privately for the F.A. examination in 1907. Though he secured a centum in 
mathematics, he failed to secure pass marks in other subjects. This marked the 
end of his formal education. \\

\medskip

\ni {\bf Formative years}\\ 

It was during the period, 1907 - 12, that Ramanujan was frantically in search of
a 
benefactor and started making contacts with those who could help him in his 
quest for a job to eke out a livelihood.  He continued to stay in Madras after his 
formal education came to an end in 1907. According to Hardy:\\

%\leftskip 1cm \rightskip 1cm 

\ni {\it The years between 18 and 25 are the critical years in a mathematician's 
 career. … During his five unfortunate years (1907-1912) his genius was   
 misdirected, side-tracked and to a certain extent  distorted}. (Hardy [1]).\\

%\leftskip 0cm \rightskip 0cm

Despite the pecuniary circumstances and the stresses and strains of day-to-day 
existence,  Ramanujan started noting down his mathematical  results in  
Notebooks.  By 1909, his Notebooks were precious to Ramanujan.  For, one 
(F.A.) classmate of his, states [13] that Ramanujan fell ill in 1909, while living in 
George Town, Madras, and on a Doctor''s advise, when he was being sent to the 
home of his parents in Kumbakonam,  Ramanujan  entrusted him with his 
Notebooks for safe keeping and stated: {\it "If I die, please hand them over to
Prof. Singaravelu Mudaliar or to the British Professor --– Edward B Ross -- Madras
Christian College}".\\

Another college mate [14] of Ramanujan has stated that during his collegiate 
years, Ramanujan taught him the method of constructing Magic Squares, the 
subject of the first chapter of his Notebooks. The interest in this subject dates 
from his school days and is disconnected from the subject matter of the remainder 
of the Notebooks. Probably Ramanujan''s expertise in preparing the conflict free 
time tables for his School inspired him to a study of these Magic Squares.\\

Ramanujan''s investigations in continued fractions and divergent series started 
during this period. His betrothal to nine year old Janaki was in 1908 and his 
wedding took place near Karur, in 1909. Robert Kanigel [8], in his biography on 
Ramanujan, constructs a vivid account of this marriage arranged  by his mother 
Komalattammal, not approved by his father, and dramatizes the foreboding of the 
impending disaster through the omens preceding the wedding, which was on the 
brink of being called off due to the late arrival of the bridegroom''s party.\\

During this period he tutored a few students in mathematics and even sought 
employment as a tutor in mathematics. Disappointed at the lack of recognition, 
during this trying period, Ramanujan had bemoaned to a friend [4] that he was 
probably destined to die in poverty like Galileo! Fortunately, this was not to
be.\\
 
In 1910,  Ramanujan sought the patronage of Prof. V. Ramaswamy Iyer --– the 
founder of Indian Mathematical Society --– who was at Salem and asked for a 
clerical job in his office. The only recommendation Ramanujan had was his 
Notebooks which by then contained several results on Magic Squares, prime 
numbers, infinite series, divergent series, Bernoulli numbers, Riemann zeta 
function, hypergeometric series, partitions, continued fractions, elliptic functions, 
modular equations, etc. A scrutiny of the entries in the Notebooks was sufficient 
to convince Prof. Ramaswamy Iyer [15] that Ramanujan was a gifted 
mathematician and he {\it "had no mind to smother his (Ramanuja'n's) genius by an 
appointment in the lowest rungs of the revenue department}". So, he sent 
Ramanujan back to Madras with a letter of introduction to Prof. P.V. Seshu Aiyar, 
then at the Presidency College, Madras. Prof. Seshu Aiyar, who had known 
Ramanujan as a student at the Government Arts College, Kumbakonam, when he
himself was employed there as a lecturer of mathematics, was meeting him after a 
gap of four years and was greatly impressed with the contents of the {\it 
well-sized} 
Notebooks. So he gave Ramanujan {\it "a note of recommendation to that true lover 
of mathematics, Dewan Bahadur R. Ramachandra Rao, who was then the District 
Collector at Nellore}."\

\eject

\ni {\bf The turning point}\\

With the help of a friend, R. Krishna Rao  [16], who was a nephew of Dewan 
Bahadur  Ramachandra Rao, Ramanujan went to Tirukkoilur in December 1910. 
This was a turning point in Ramanujan's life. Ramachandra Rao states [17] that  
{\it "in the plentitude of my mathematical wisdom, I condescended to permit 
Ramanujan to walk into my presence}". At that time, Ramanujan appeared to  
Ramachandra Rao as \\

%\leftskip 1cm \rightskip 1cm

\ni {\it a short uncouth figure, stout, unshaved, not over-clean, with 
one conspicuous feature - shining eyes - walked in, with a frayed Notebook under 
his arm $\cdots$. He was miserably poor. He had run away from Kumbakonam to get 
leisure in Madras to pursue his studies. He never craved for any distinction. He 
wanted leisure, in other words, simple food to be provided for him without 
exertion on his part and that he should be allowed to dream on.}\\

%\leftskip 0cm \rightskip 0cm

 Though Ramachandra Rao gave him a patient hearing, he took a few days to look 
into the Notebooks of Ramanujan. At their fourth meeting, when Ramanujan 
confronted Ramachandra Rao with a letter from Prof. Saldhana of Bombay 
appreciating the genuineness of his work, Ramachandra Rao started to feel that 
Ramanujan's work must be examined in depth by eminent mathematicians. 
Ramachandra Rao himself states [17] that Ramanujan led him {\it step-by-step to 
elliptic integrals and hypergeometric series and at last to his theory of divergent 
series not yet announced to the world} and this converted him into a benefactor 
who undertook to underwrite Ramanujan's expenses at Madras for some time. \\

Prof. Seshu Aiyar also communicated the earliest contributions of Ramanujan to 
the Journal of the Indian Mathematical Society (I.M.S.) in the form of questions. 
These appeared in 1911 and in his brief and illustrious career Ramanujan 
proposed in all 59 questions or solutions to questions in this journal. The first
fifteen page article entitled: {\it Some properties of Bernoulli numbers} appeared
in the 
same 1911 volume of the journal of the I.M.S. In it Ramanujan stated eight 
theorems embodying arithmetical properties of the Bernoulli numbers, indicating
proofs for three of them; two theorems are stated as corollaries of two others, 
while three theorems are stated as mere conjectures. Prof. Seshu Iyer states [18]: 
{\it Ramanujan's methods were so terse and novel and his  presentation was 
so lacking in clearness and precision, that  the ordinary reader, unaccustomed to 
such intellectual gymnastics, could hardly follow him.}\\

Ramanujan lived in a small house, called `'Summer Hous'e', in Sami Pillai Street, 
Triplicane, Madras, accepting reluctantly a monthly financial assistance from the 
collector of Nellore for about a year. Later he declined this help and from Jan. 12 
to Feb. 21, 1912, he worked as a clerk in the Accountant General's Office, on
a salary of Rs.25/- per month.  Not  satisfied with this job, Ramanujan applied for 
and secured a post in the Accounts Section (Class III, Grade IV clerk on a salray 
of Rs.30/- per month) in the Madras Port Trust, with the help of Mr. S. Narayana 
Iyer, the Manager of Port Trust, who was the treasurer of the IMS and a friend of 
Profs. V. Ramaswamy Aiyar and  P.V. Seshu Aiyar.\\

 Mr. Narayana Aiyer was a good mathematician and was a great source of support 
to Ramanujan. He was not only instrumental in Ramanujan  being offered a job in 
the Madras Port Trust, but also in securing  for Ramanujan the life-long support
of Sir Francis Spring. When Ramanujan was living in No. 580, Pycrofts Road, 
Triplicane, Madras, he used to meet Mr. Narayana Iyer and work out Mathematics 
on two big slates. Narayana Aiyer's son N. Subbanarayanan relates the role his 
father played in the career  of Ramanujan [18, p. 112]:\\

%\leftskip 1cm \rightskip 1cm

\ni {\it My father, being a fairly good mathematician himself, was unable 
        to capture the strides of Ramanujan's discoveries. He used to tell him,      
       ``"When I am not able to understand your steps, I do not know how other        
       mathematicians of a critical nature will accept your genius. You must   
       descend to my level and write at least ten steps between the two steps of   
       yours''.'' Sri Ramanujan ud
       to say, ``"When it is so simple and clearto me, 
       why should I write more steps ?'''' But somehow my father slowly got him     
       round, cajoled him and made him write some more, though it used to be a 
       mighty task of boredom to him.}\\

%\leftskip 0cm \rightskip 0cm

Dewan Bahadur Ramachandra Rao wrote to Sir Francis Spring, Chairman of 
Madras Port Trust, about Ramanujan. He also induced Prof. C.L.T. Griffith of
the Engineering College, Madras to take interest in Ramanujan and Prof. Griffith 
in turn wrote [19] in  November 1912, to Sir Francis Spring, the Chairman of 
Madras Port Trust about the  very poor accountant who was a most remarkable 
mathematician and asking him to keep Ramanujan  happily employed until 
something can be done to make use of his extraordinary gifts. As stated before, 
these efforts resulted in Ramanujan''s entry into Port Trust, on March 1, 1912,  as 
a Clerk in the Accounts Department. This may well be considered as 
{\underline{the}} turning 
point in his career prospects. He held this clerical post for 14 months. His wife 
joined him during this period and Ramanujan shifted his residence to Saiva 
Muthiah Mudali Street in George Town. This period also marked the beginning of 
the appreciation of his scholarship and researches in mathematics.  \\

Prof. Griffith  wrote to Prof. M.J.M. Hill, of University College, University of 
London, on Ramanujan''s work and he received a reply in December 1912. 
Unfortunately, Prof. Hill [20] could not find time to study the results. He observed 
that {\it 'the book which will be most useful to him is Bromwic'h's Theory of
Infinite 
Series, published by Cambridge University Press (or Macmillan)}' and gave advice 
as to how Ramanujan could get his papers published. In a sequel to this reply, 
dated 7 December 1912, Prof. Hill wrote to Prof. Griffith [21]: 

%\leftskip 1cm \rightskip 1cm

\ni {\it Mr. Ramanujan is evidently a man with a taste for Mathematics,  and with
some 
ability, but he has got on the wrong lines. He does not understand the precautions 
which have to be taken in dealing with divergent series, otherwise he could not 
have obtained the erroneous results you send me, viz.

\vskip -7mm

\begin{eqnarray*} 
1\ + 2\ + 3\ +\ \cdots     +\ \infty     & = &  -1/12,\\
1^2 + 2^2 + 3^2 + \cdots  …  +\infty^2 &=& 0, \\
1^3 + 2^3 + 3^3 + \cdots … +\infty^3   &=& 1/240.
\end{eqnarray*}

\ni The sums of  $n$  terms of these series are:
$$ n(n+1)/2,\qquad  n(n+ 1/2 )(n+1)/3,\qquad       [n(n+1)]^2 / 2 $$
and they all tend to $\infty$  as  $n$  tends to $\infty$ . I do think you can do
no better for him 
than to get him a copy of the book I recommended, Bromwich's Theory of Infinite 
Series, published by Macmillan and Co., who have branches in Calcutta and 
Bombay. Price 15/- net.}\\

%\leftskip 0cm \rightskip 0cm

It is not as though Ramanujan was not aware of the apparent absurd looking 
nature of the results on divergent series.  Ramanujan, in his second letter  to  
Hardy [22], wrote: \\

%\leftskip 1cm \rightskip 1cm

\ni {\it I have got theorems on divergent series, theorems to calculate the convergent 
values corresponding to the divergent series, viz.:

\vskip -7mm

\begin{eqnarray*}
1 - 2 + 3 - 4   +\cdots & = &    1/4, \\
1 -   1! + 2! -    3! + \cdots & = &  0.596,\\
1 + 2 + 3 + 4  + \infty & = & -1/12,\\    
1^3 + 2^3 + 3^3 + \cdots  +\infty^3& = &  1/24.
\end{eqnarray*}

\ni Theorems to calculate such values for any given series (say, $1 - 1 ^1 + 2^2 -
3^3+ 4^4 - 5^5 + \cdots $), 
and the meaning of such values. I have also dealt with such questions 
'When to use, where to use, and how to use such values, where do they fail and 
where do they not ?'}\\

%\leftskip 0cm \rightskip 0cm 

\ni Hill failed [23] to discern the origin of the results of Ramanujan and the
three 
sums of the integers, their squares and their cubes are indeed the values of
$\zeta(-n)$, for $n = 1,2,3$, respectively\footnote{The 
Riemann zeta function is defined as: $\zeta(s) = \sum_{n=1}^\infty  1/n^s, Re\
s > 1$. The $\zeta$ 
function has a  unique analytic continuation to the points $s=-1$,
where we get $\zeta(-1) = -1/12$,
which is what Ramanujan writes as: $1+2+3+\cdots + \infty = -\frac{1}{12}$. This
result is used in the zeta
function regularization method, by String theorists, in recent times.}.\\  

Ramanujan
published two short notes, one {\it On question 330 of Professor Sanjana} and
another a
{\it Note on a set of simultaneous equations}, in the IMS journal, in 1912. When
Ramanujan
approached Prof. Seshu Aiyar with some theorems on Prime Numbers, his attention was
drawn to G.H. Hardy's Tract on{\it  Orders of infinity}. In it, Ramanujan observed
that
([III], p.xxii): {\it "no definite expression has yet been found for the number of
prime
numbers less than any given number}". Ramanujan told Prof. Seshu Aiyar that he ha
discovered the required result. This made Prof. Seshu Aiyar suggest communication
of this and other results to Mr. G.H. Hardy -- a Fellow of the Royal Society and
Cayley Lecturer in Mathematics at Cambridge – a world famous mathematician, who was
ten years Ramanujan''s senior. \\

\medskip

\ni {\bf The years of fruition}\\

The life of Ramanujan, in the words of C.P. Snow [24] {\it is an admirable story,
and 
one which showers credit on nearly everyone} [25]. Ramanujan''s first letter [26] to 
Prof. Hardy, dated 16th January 1913, is a historic letter. It contained the bare 
statements of about 120 theorems, mostly formal identities from the Notebooks. 
This collection obviously represented what Ramanujan himself considered were 
results of importance.  Ramanujan wrote:\\

%\leftskip 1cm \rightskip 1cm

\ni {\it Dear Sir,\\

     I beg to introduce myself to you as a clerk in the Accounts Department of the 
Port Trust Office at Madras on a salary of  \pounds 20 per annum. I am now about 23 
years of age. I have had no University education but I have undergone the 
ordinary school course. After leaving school I have been employing the spare 
time at my disposal to work at Mathematics. I have not trodden through the 
conventional regular course which is followed in a University course, but I am 
striking  out a new path for myself. I have made a special investigation of 
divergent series in general and the results I get are termed by the local 
mathematicians as `startling' $\cdots$.\\

    I would request you to go through the enclosed papers. Being poor, if you are 
convinced that there is anything of value I would like to have my theorems 
published. I have not given the actual investigations nor the expressions that I get 
but I have indicated the lines on which I proceed. Being inexperienced I would 
very highly value any advice you may give me. Requesting to be excused for the 
trouble I give you,\\
  
    I remain,    Dear Sir,
                                                     
\hfill Yours truly,

\hfill                     (sd) S. Ramanujan}.\\

%\leftskip 0cm \rightskip 0cm

Prof. Hardy, the professional mathematician, who was aware that he was  
{\it the first 
really competent person who had the chance to see some of his work}, found some 
of the series formulae intriguing, some of the integral formulae (which were 
classical and known) vaguely familiar and he could prove some integral formulae 
with effort but these were to him the least impressive.  However, some of 
Ramanujan''s formulae were{\it on a different level} and obviously both difficul
and deep, which even Hardy [27] {\it  'had never seen anything in the least like them 
before}' and whic he has state {\it `defeated me completely'}.\\

The following is a record of  Hardy''s reaction to this historic letter of
Ramanujan, in the words of C.P. Snow [28] : \\

%\leftskip 1cm \rightskip 1cm

\ni {\it Hardy gave the manuscript a perfunctory glance, and went on reading the 
morning paper. It occurred to him that the first page was a little out of the 
ordinary for a cranky correspondent. It seemed to consist of some theorems, very 
strange-looking theorems, without any argument. Hardy then decided that the 
man must be a fraud, and duly went about the day according to his habits, giving 
a lecture, playing a game of tennis. But there was something nagging at the back 
of his mind. Anyone who could fake such theorems, right or wrong must be a 
fraud of genius. Was it more or less likely that there should be a fraud of genius 
or an unknown Indian mathematician of genius ? He went that evening after 
dinner to argue it out with his collaborator, J.E. Littlewood, whom Hardy always 
insisted was a better mathematician than himself. They soon had no doubt of the 
answer. Hardy was seeing the work of someone whom,  for natural genius, he 
could not touch – who, in natural genius, though of course not in achievement, as 
Hardy said later, belonged to the class of Euler and Gauss.\\

\ni Hardy made up his mind that Ramanujan should be brought to Cambridge and 
provided with the necessary education and contact with western mathematicians 
of the highest class.  So, Hardy,  wrote to the Secretary of the Indian students, in 
the India Office, London, suggesting that some means be found to get Ramanujan 
to Cambridge and he in turn wrote, in February 1913, to Mr. Arthur Davies, the 
Secretary to the Advisory Committee for Indian students in Madras conveying the 
desire of the tutors at Trinity to get Ramanujan to Cambridge.}\\

%\leftskip 0cm \rightskip 0cm

Sir Francis Spring, the Chairman and Mr. S. Narayana Iyer, the Manager of 
Madras Port Trust gave Ramanujan every possible encouragement. Dr. Gilbert T. 
Walker, F.R.S., Director General of Observatories, Simla,  and Head of the Indian 
Meteorological Department, paid a visit to the harbour in Madras on February 25,
1913 and Sir Francis Spring drew his attention to the work of  Ramanujan and his 
Notebooks. Dr. Walker, a good mathematician and a Senior Wrangler, was a 
former Fellow of Trinity College, Cambridge, as well as a lecturer and he said 
that in his opinion Mr. Hardy would be the most competent to arrive at a 
judgement of the true value of the work of Ramanujan. Since by then Hardy''s 
reply had arrived (on Feb. 8, 1913),  Gilbert Walker  wrote [29] to Mr. Francis 
Dewsbury, the Registrar of  the University of Madras, commending the work of 
Ramanujan to be {\it comparable in originality with that of a Mathematics Fellow in 
a Cambridge college}, though lacking in the precision and completeness necessary 
for establishing the universal validity of the results. He wrote that it was perfectly 
clear to him that {\it the university would be justified in enabling S. Ramanujan
for a 
few years at least to spend the whole of his time on mathematics without any 
anxiety as to his livelihood}. He also wanted the University to correspond with Mr. 
Hardy, Fellow of Trinity College, Cambridge, since Ramanujan was already in 
correspondence with Hardy, assuring Mr. Hardy of the University's interest in 
Ramanujan. The recommendation of Dr. Walker was accepted by the Board of 
Studies in Mathematics of the University of Madras. Then the Vice Chancellor of 
the University got the approval  of the Syndicate overcoming the legal hurdle of 
awarding a research scholarship to  Ramanujan who did not have the required 
qualification of a Master's Degree. As a measure of  precaution, the consent of the 
Chancellor of the University (Lord Pentland, the Governor of Madras) was 
obtained to grant Ramanujan   a special research scholarship of Rs.75/- per month 
for two years with the condition that Ramanujan should submit quarterly reports 
on his work . The Madras Port Trust granted Ramanujan two years leave (on loss 
of pay) to enable him to accept this scholarship from May 1913, as the {\bf first 
research scholar of the University of Madras}. Thus began Ramanujan''s carer
as a professional mathematician.\\

In quick succession, Ramanujan received in the next three months, four long 
letters [30] from Hardy in which the latter wrote plainly about what had been 
proved or claimed to have been proved by Ramanujan. He clearly communicated 
his genuine anxiety {\it to see what can be done to give you (Ramanujan) a better 
chance of making the best use of your obvious mathematical gifts}. At last 
Ramanujan had found a sympathetic friend in Hardy and was willing to place 
unreservedly in his hands all that he had.\\

Ramanujan  wrote again to Hardy on 27th February 1913 and sent him more 
formulae and explanations.  On 17th April 1913, Ramanujan wrote to Hardy 
about his having secured the scholarship, of \pounds 60 per annum, of the
University of 
Madras, for two years. Ramanujan took up residence at Hanumantharayan Koil 
Lane in Triplicane around this time and had  access to books on mathematics in 
the University library. His wife Janaki and his mother came to live with him.\\

Ramanujan was initially reluctant to go abroad because of his own caste 
prejudices\footnote{Crossing the oceans was considered a sacrilege by the Hindu
Brahmins and often
people did so were, on their return to India, treated as outcastes. All
relationships with the even their families were shunned!} 
in those days which were compounded by the extremely orthodox views of his
mother to whom he was greatly attached. At the beginning of 1914, Mr. E.H. Neville,
a young mathematician and a Fellow of Trinity College, Cambridge, was in Madras as
a visiting lecturer to give a series of lectures on Differential Geometry to
Mathematics Honours students of the University of Madras. Mr. Hardy entrusted him
with the mission of persuading Ramanujan to visit Cambridge. Mr. Neville met
Ramanujan and saw his priceless notebooks. This was sufficient to convince him of
Ramanujan''s uncommon ability and to make him take over the initiative to overcome
all the difficulties in arranging for Ramanujan''s visit to Cambridge. Prof.
Richard
Littlehails, who was a Professor of Mathematics with the observatory in Madras
introduced Neville [31] to everyone who carried weight in the University or in the
civil administration. Neville, in turn, explained to them the importance of
Ramanujan''s stay in Cambridge, and urged them to be generous in their support.\\

In a letter [32], dated 28th January 1914, to Mr. Dewsbury, the Registrar of the 
University of Madras, Mr. Neville wrote about {\it the importance of securing to 
Ramanujan a training in the refinements of modern methods and a contact with 
men who know what range of ideas have been explored and what have not} and 
prophesied that Ramanujan would respond to such a stimulus and that {\it his name 
will become one of the greatest in the history of mathematics, and the University 
and city of Madras will be proud to have assisted in his passage from obscurity to 
fame}. The very next day, Prof. Littlehails also wrote [33] to Mr. Dewsbury that 
Ramanujan {\it  be granted by this University a scholarship of about \pounds 250
(Sterling) 
together with a grant of about \pounds 100 in order to enable him to proceed to 
Cambridge. Ramanujan is a man of most remarkable mathematical ability, 
amounting I might say to genius, whose light is metaphorically hidden under a 
bushel in Madras}.\\ 

The proposals regarding the scholarship to be granted to Ramanujan by the 
University of Madras were approved. To the lasting credit of the University of 
Madras, the Syndicate decided within a week to set aside Rs.10,000/- to offer 
Ramanujan a scholarship of \pounds 250 a year plus \pounds 100 for a passage by
ship and for initial outfit\footnote{The second class fare between London and
Bombay was \pounds 32, in 1914, or about Rs. 480 -- {\it British Passenger Liners
of the Five Oceans}, By C.R. Vernon Gibbs (London: Putnam, 1963, p.63). ([8],
p.397).}.
At the instance of Professors Neville and Littlehails, Sir Francis 
Spring wrote [34] to the personal Secretary (Mr. C.B. Cotterell) to the Governor 
(Lord Pentland) of Madras, persuading His Excellency to speedily approve the 
University's sanction.  Government sanction too was granted within a wek.\\

This offer of the University of Madras was made to Ramanujan in February 1914. 
He sent his wife and mother back to Kumbakonam, changed the traditional hair-
style of  a  brahmin, viz. a tuft, and  got his hair trimmed in European style and 
left Madras by s.s. Nevasa on 17th March 1914. Prior to his departure, he 
arranged with the University for \pounds 60 a year to be sent to his parents in 
Kumbakonam, out of his annual scholarship amount. Mr. Arthur Davies and Prof. 
Littlehails attended to all the details regarding Ramanujan''s passage to England. 
Except for the first three days when he was sea-sick, Ramanujan enjoyed the 
voyage and reached London through the Channel and the Thames on 14th April 
1914. He was received by Mr. E.H. Neville and his brother at the docks and 
stayed at Cromwell Road for a few days before going to Cambridge on the 18th 
evening. He remained for a few days in Mr. Neville''s house before moving to th
college premises for stay, which even though costlier than lodging houses, was 
more convenient for him and the professors. Ramanujan wrote [35] to his friend 
that {\it Mr. Hardy, Mr. Neville and others here are unassuming, kind and obliging. 
As soon as I came here, Mr. Hardy paid \pounds 20 to the college for my entrance
and 
other fees and made arrangements to give me a scholarship of \pounds 40 a year}.\\

Ramanujan was admitted by Mr. Hardy to Trinity College which supplemented 
his scholarship with the award of an exhibition of \pounds 60 a year, to augment
the \pounds 250 a year scholarship  awarded by the University of Madras.\\

Though Ramanujan had access only to Carr's Synopsiss--– and perhaps, to a few 
other books\footnote{From the article of Mr. Narayana Iyer's son [Ref. 18, p.112],
Ramanujan had access to
a book on Jacobis elliptic functions. Unfortunately, it is not possible to
ascertain, from the records of the Library of the University of Madras, what books
were available for reference to Ramanujan.} -- still, in the words of the
historian J.R. Newman [36], he { \it 
  arrived in England abreast and often ahead of contemporary mathematical
knowledge.  Thus, in a lone mighty sweep, he had succeeded in recreating in his
field, through his own unaided powers, a rich half century of European mathematics.
One may doubt whether so prodigious a feat had ever before been accomplished in the
history of thought. }\\

To Mr. Hardy [37] Ramanujan''s friend, philosopher an discoverer: \\

%\leftskip 1cm \rightskip 1cm 

\ni {\it  The
limitation of his knowledge was as startling as its profundity. Here was a man
who could  work out modular equations, and theorems of complex multiplications, 
to orders unheard of, whose mastery of continued fractions was, on the formal 
side at any rate, beyond that of any mathematician in the world, who had found 
for himself the functional equation of the zeta-function, and the dominant terms of 
many of the most famous problems in the analytic theory of numbers, and he had 
never heard of a doubly periodic function or of Cauchy's theorem, and had 
indeed but the vaguest idea of what a function of a complex variable was. His 
ideas of what constituted a mathematical proof were of the most shadowy 
description. All his results, new or old, right or wrong, had been arrived at by a 
process of mingled argument, intuition and induction, of which he was entirely 
unable to give a coherent account.}\\

%\leftskip 0cm \rightskip 0cm

With such a natural genius, Hardy collaborated and tried to teach, as he wrote, the  
{\it things of which it was impossible that he should  remain in ignorance. … It
was 
impossible to allow him to go through life supposing that all the zeroes of the zeta 
function were real. So I had to try to teach him, and in a measure I succeeded, 
though I obviously learnt from him much more than he learnt from me} [38].\\ 

Hardy did not attempt to convert Ramanujan into a mathematician of the modern 
school but enabled him to go on producing  original ideas in his classical mould 
with rigorous proofs for the theorems he discovered.\\

The period of Ramanujan''s stay in England almost overlapped with the years in
which World War I took place.  {\it One of the lecturers went to 
war\footnote{Ramanujan was perhaps referring to the departure of Mr. J.E.
Littlewood.}} 
wrote Ramanujan [39] to a friend in India and Ramanujan felt that {\it the other
professors $\cdots$ lost their interest owing to the $\cdots$  war}. One of the
professors had 
remarked that Ramanujan was in England at the most unfortunate time. There 
were about 700 students before the war, but this number was reduced to 150 by 
November 1915.\\

Initially Ramanujan asked for and obtained some South Indian food items (like 
tamarind, coconut oil, etc.) by post parcel from his home, as well as from a 
company in London but by January 1915, he wrote [23] to a friend of his in India 
that {\it  now as well as in the future I am not in need of anything as I gained
control 
over my taste and can live on mere rice with a little salt and lemon juice for an 
indefinite time.}  His difficulty of getting proper food was alleviated by the 
availability of good milk and fruits. Being a vegetarian he had no option but to 
cook for himself.\\

He was attending a lecture by Mr. Berry at the University on elliptic integrals. Mr. 
Berry was working out some formulae on the black-board and a glance at 
Ramanujan''s face, alight with excitement, caused him to ask Ramanujan whether 
he was following the lecture and whether he had anything to say. At this
Ramanujan went to the black-board and much to everyone's surprise wrote down 
some of the results which were yet to be proved. This anecdote was recalled by 
Dr. P.C. Mahalanobis [40], the eminent Indian statistician, who joined King''s 
College, Cambridge, in October 1913, and took a mathematics course by Prof. 
Hardy.  The following is another anecdote about Ramanujan from Dr. 
Mahalanobis [40]: {\it I was fortunate in forming a good friendship with Ramanujan 
very  soon.  It came about in a somewhat strange way. One day, soon after his 
arrival, I went to see Ramanujan in his room in Trinity College. It had turned 
quite cold. Ramanujan was sitting near the fire\footnote{Ramanujan's room had
electricity and he was provided with a gas stove.}. I asked him whether he was
quite 
warm at night. He said that he was feeling the cold though he was sleeping with 
his overcoat on and was also wrapping himself up in a shawl. I went to his 
bedroom to see whether he had enough blankets. I found that his bed had a 
number of blankets but all tucked in tightly, with a bed cover spread over them. 
He did not know that he should turn back the blankets and get into the bed. The 
bed cover was loose; he was sleeping under that linen cover with his overcoat  
and shawl. I showed him how to get under the blankets. He was extremely 
touched. I believe this was the reason why he was so kind to me.}\\

Ramanujan wrote a few articles soon after he reached Cambridge and in June 
1914, Hardy presented some of the results from Ramanujan''s Notebooks at a 
meeting of the London Mathematical Society. However, in January 1915, 
Ramanujan wrote [41] to a friend in India that his notebook is sleeping in a corner 
for these four or five months. Ramanujan was more interested in getting 
new results (and partly due to the ongoing war),  he decided to publish the old 
results worked out in his Notebooks after the war. After about a year and a half at 
Cambridge, Hardy wrote to the Registrar of the University of Madras, that 
Ramanujan {\it  is beyond question the best Indian mathematician of modern times.  
He will always be rather eccentric in his choice of subjects and methods of 
dealing with them. But of his extraordinary gifts there can be no questions; in 
some ways he is the most remarkable mathematician I have ever known. }\\

Hardy''s letter [42] and official report to the University, as well as an appeal by
Sir 
Francis Spring to the University to continue the assistance extended by it to 
Ramanujan, made the University (in December 1915) extend the scholarship up to 
March 1919.\\

\medskip

\ni {\bf Honours} \\

During his five year stay in Cambridge, Ramanujan published 21 research papers 
containing theorems on definite integrals, modular equations, Riemann''s zeta 
function, infinite series, summation of series, analytic number theory, asymptotic
formulae, modular functions, partitions and combinatorial analysis. His paper 
entitled  {\it Highly Composite Numbers} which appeared in the Journal of the
London 
Mathematical Society, in 1915, is 62 pages long and contains 269 equations. This 
is his longest paper. The London Mathematical Society had some financial
difficulties at that time and Ramanujan was  requested to reduce the length of his 
paper to save printing expenses. Five of these 21 research papers were in 
collaboration with Hardy. Ramanujan also published  5 short notes in the Records 
of Proceedings at meetings of the London Mathematical Society and six more in 
the journal of the Indian Mathematical Society.\\

Ramanujan was awarded the B.A. degree by research in March 1916 for his work 
on Highly composite numbers and published as a long paper. Ramanujan''s
dissertation bore the same title and included six other papers. Ramanujan was 
registered as a research student in June 1914 and the prerequisite of a diploma or 
a certificate, as well as the domiciliary requirement of six terms must have been 
relaxed in his extraordinary case. It is unfortunate that a copy of this dissertation 
is not to be found in the records of the University [43]. According to Hardy [44], 
this work of Ramanujan {\it is a very peculiar one, standing somewhat apart from 
the main channels of mathematical research. But there can be no question as to 
the extraordinary  insight and ingenuity which he has shown in treating it, nor 
any doubt that the memoir is one of the most remarkable published in England for 
many years.}\\

Ramanujan''s designated tutor who monitored his progress at Trinity College, 
Cambridge, was E.W. Barnes, who considered Ramanujan as perhaps the most 
brilliant of all the top Trinity students, which included Littlewood [45]. Hardy 
was immensely satisfied with the progress of Ramanujan and wrote so to the 
Registrar of the University of Madras supporting an extension of Ramanujan''s 
two-year scholarship {\it $\cdots$ "… until, as I confidently expect, he is elected
to a 
Fellowship at the College. Such an election I should expect in October 1917".}
Later, in June 1916, in an official report on the progress of  Ramanujan's work in 
England to the University's Registrar, he wrote: {\it $cdots$ " …it is already safe
to say that 
Mr. Ramanujan has justified abundantly all the hopes that were based upon his 
work in India, and has shown that he possesses powers as remarkable in their 
way as those of any living mathematicians. $\cdots$  I have said enough, I hope,
to give some idea of his astonishing individuality and power. India has produced many 
talented mathematicians in recent years, a number of whom have come to 
Cambridge and attained high academical distinction. They will  be the first to 
recognize that Mr. Ramanujan's work is of a different category".}\\

In spite of the war which was raging, which deprived Ramanujan of  the center 
stage which he would otherwise have held with his brilliant research work in the 
midst of his peers, the confidence he kindled in Hardy was enough to win for him 
recognition and laurels very soon, but, unfortunately, the first signs of illness 
appeared in Ramanujan in the spring of 1917.\\

Thanks to the unstinted efforts of Hardy, who did his best to get Ramanujan due 
recognition, he was elected a Fellow of the Royal Society of London in February 
1918. The Records of the Royal Society, dated December 18, 1917, include the 
following certificate for the candidature of Ramanujan (then  a Research student 
in Mathematics at Trinity College, Cambridge) for election to the Fellowship of 
the Royal Society\footnote{A copy of this documjent is an exhibit in the Ramanujan
Museum in Royapuram, Madras.}:\\

\ni {\it Qualifications (Not to exceed 250 words):\\

\ni  Distinguished as pure mathematician, particularly for his investigations in
elliptic functions and the theory of numbers. Author of the following papers, 
amongst others: `Modular equations and approximations to $\pi$',
{\underline{Quarterly 
Journal}}, vol. 45; `New expressions of Riemann's function $\zeta(s)$ and $\chi
(t)$'',ibid, 
vol. 46; `Highly composite numbers', {\underline{Proc. London. Math. Soc.}}, vol.
14; `On 
certain arithmetical functions', {\underline{Trans. Camb. Phil. Soc.}}, vol. 22;
`On the 
expression of a number in  the form  a $x^2 +be y^2 +c z^2 +d t^2$,
{\underline{Proc. Camb. Phil. Soc.}} , vol. 19.\\

\ni Joint author with G.H. Hardy, F.R.S., of the following papers: `Une formulae 
asymptotique pour le nombre des partitions de n',
{\underline{Comptes Rendus}} , 2 Jan. 1917; 
`Asymptotic Formulae for the distribution of numbers of various types',
{\underline{Proc. London Math. Soc.}}, 
vol. 16; `The normal number of prime factors  of a number 
n', {\underline{Quarterly Journal}}, vol. 47; `Asymptotic Formulae in Combinatory
Analysis',{\underline{Proc. London Math. Soc.}}, (awaiting publication).\\

\ni being desirous of admission into the ROYAL SOCIETY OF LONDON, we the
undersigned propose and recommend him as deserving that honour, and as likely
to become a useful and valuable Member.}\\

This nomination was  proposed by  G.H. Hardy and seconded by  P.A. 
MacMahon. The signatories with `Personal knowledge'' of Ramanujan were, 
besides Hardy and MacMahon,   J.H. Grace, Joseph Larmor, T.J.I''A. Bromwich   
E.W. Hobson\footnote{Note that E.W. Hobson and H.F. Baker, who had not
replied to letters written by Ramanujan from India, being signatories.}, H.F.
Baker, J.E.Littlewood and J.W. Nicholson. Besides these 9 
signatures were the signatures of E.T. Whittaker, A.R. Forsyth and A.N. 
Whitehead, under those who knew him  from General Knowledge. This 
certificate on a printed form of the Royal Society has been filled by hand (and the 
hand writing appears to be that of Mr. Hardy), delivered at the Apartments of the 
Society on the 18th  Dec. 1917 and read to the Society on the 24th January 1918.\\

As a consequence, Ramanujan was, awarded on Feb. 28, 1918, the Fellowship of 
Royal Society, London, and the citation read:\\

%\leftskip 1cm \rightskip 1cm

\ni {\bf Srinivasa Ramanujan, Trinity College, Cambridge. Research student in     
Mathematics Distinguished as a pure mathematician particularly for  
his investigations in elliptic functions and the theory of numbers.}\\

%\leftskip 0cm \rightskip 0cm

In recent times\footnote{Private communication by e-mail from Prof. Dalitz,
Oxford University, dated March 29, 1996.}, 
I came to know from Prof. R.H. Dalitz, F.R.S., that the signature 
of Ramanujan is not in the book of the Royal Society. According to Prof. Dalitz: 
{\it The book is indexed, so it is just not there. The reason undoubtedly is that
he was 
ill in that period and could not go to the Royal Society to sign it. There are other 
examples of well-known F.R.S's who somehow didn't get their signature into the 
book. That means that he did not ever attend any meeting of the Royal Society;  if 
he had, they would have brought out the book and not let him go until he had 
signed. Of course, it was also war-time, which meant that there were as few 
meetings as possible.}\\

Ramanujan was elected to a Trinity College Fellowship, in October 1918, which 
was a Prize Fellowship worth \pounds 250 a year for six years with no duties or 
conditions.  These awards acted as great incentives to Ramanujan who discovered 
some of the most beautiful theorems in mathematics, subsequently.\\
 
Hardy's letter to the Registrar of the University of Madras, Mr. Dewsbury,
dated Nov. 26, 1918 [46] struck a hopeful note: {\it There is at last, I am
profoundly 
glad to say, a quite definite change for the better. I think we may now hope that he 
has turned the corner, and is on the road to recovery. His temperature has ceased 
to be irregular, and he has gained nearly a stone\footnote{One stone weight 
is equal to 14 pounds.} in weight. The consensus of 
medical opinion is that he has been suffering from some obscure source of blood 
poisoning, which has now dried up; and that it is reasonable to expect him to 
recover his health completely and if all goes well fairly rapidly.}\\

Ramanujan''s symptoms were predominantly night-time fever, loss of weight 
leading to his emaciated looks and these caused depressions which once drove 
him to the limit of attempting suicide\footnote{A story which was recounted many 
years after his death, by the Astrophysicist Dr. 
S. Chandrasekhar, Nobel Laureate, as told to him by Prof. Hardy, and reproduced in
Ch. 5 of Ref. 7.}. These symptoms made the doctors consider various diagnosis, at
different times: 
gastric ulcer, malaria, tuberculosis, cancer of the liver, etc.  In recent times,
with hind sight, vitamin $B_{12}$ deficiency (something unknown to the world at
that
time) has been diagnosed as a possibility [47]. The recovery alluded to by Hardy in
his letter to Dewsbury was obviously the reason why Ramanujan was persuaded to
return to India, with the hope that he would soon recover and return to take up the
Trinity College Fellowship awarded to him for five years.\\

\medskip

\ni {\bf The beginning of the end}\\

After completing nearly five years at Cambridge, early in 1919, when Ramanujan 
appeared to have recovered sufficiently to withstand the rigours of a long voyage 
to India, he left England on 27th February 1919 by s.s. Nagoya. Four weeks later 
on 27th March he arrived at Bombay and soon after at Madras, thin, pale and 
emaciated, but {\it with a scientific standing and reputation such as no Indian has 
enjoyed before}. Professor Hardy who expressed this view [48] also hoped that 
{\it India will regard him as the treasure he is}. He urged the University of
Madras to 
make a permanent provision for him to enable him to continue his research work. 
Again the University rose to the occasion by granting Ramanujan \pounds 250 a year
as 
an allowance for five years, commencing from April 1919. He was sent back to 
India by Hardy with the fond hope that the warmer climate would help complete 
his recovery from a tubercular tendency.\\

Most unfortunately his precarious health did not improve, on his return to India. 
Fevers relapsed and in addition, his wife recalled that he suffered severe bouts of 
stomach pain too [49]. Ramanujan was subject to fits of depression, had a 
premonition of his death and was a difficult patient. He spent 3 months in Madras, 
2 months in Kodumudi and 4 months at Kumbakonam. When his condition 
showed signs of further deterioration, after great persuasion, Ramanujan was 
brought to Madras for expert medical treatment, in January 1920. Despite all the 
tender attention he could get from his wife who nursed him throughout this 
period, and the best medical attention from the doctors, his  untimely end came on 
26th April 1920, at Chetput, Madras, when Ramanujan was 32 years, 4 months 
and 4 days old.  His wife lived with him, after she came of age, only for a year 
before his departure to England, and looked after him during his illness after his 
return. Even during those months of prolonged illness he kept on working, though 
in a reclining position, at a furious pace and kept jotting down his results on 
sheets of paper. In his last and only letter to Hardy written after his return to India, 
in January 1920, Ramanujan communicated his original work on what he called 
`'mock'' theta functions. \\

From the available evidence and retrospective diagnosis, Young [59] makes out
the 
case for ``'hepatic amoebiasis''', a tropical disease contacted by Ramanujan in
1906, 
as the cause of his terminal illness. His reason as to why this was not recognized 
at that time is best  recounted in his own words: {\it Hepatic amoebiasis was
regarded 
in 1918 as a tropical disease (`'tropical lie r abscess''), and this would have
had 
important implications for successful  diagnosis, especially in provincial medical 
centers. Furthermore, the specialists called in were experts in either tuberculosis 
or gastric medicine. Another major  difficulty is that a patient with this disease 
would not, unless specifically asked, recall as relevant that he had had two 
episodes of dysentery 11 and 8 years before. Finally, there is the very good 
reason that, because of the great variability in physical findings, the diagnosis 
was difficult in 1918 and remains so today:  hepatic amoebiasis 'presents a
severe challenge to the diagnostic skills … [and] should be considered in any 
patient with fever and an abnormal abdominal examination coming from an 
endemic area'.}\\
  
Hardy, who was unaware that the end was to come so soon was shocked when it 
came prematurely. He was of the view that a mathematician is often 
comparatively old at 30.  For, in his roll-call of mathematicians, Hardy wrote 
([50], p.71):  {\it Galois died at twenty-one, Abel at twenty-seven, Ramanujan at 
thirty-three, Riemann at forty … I do not know an instance of a major
mathematical advance initiated by a man past  fifty\footnote{A few examples which
can be cited which explode `The Myth of the Young
Mathematician'
are: Newton's Principia was written when he was in his mid 40s;  when Euler, despite
his blindness, produced his three volumes on integral calculus when 
he was in his 60s; Gauss at 34 proposed his theory of analytic functions; and in
more recent times, Cartan, Poincar\'e , Siegel, Kolmogorov and Erd\"os exhibited 
creativity in mathematics in their later years. (Ref. Susan Landau, Notices
of the AMS, vol. 44 (1997) p. 1284.)}.}\\

\medskip

\ni {\bf Human qualities}\\

{\it In figure he (Ramanujan) was a little below medium height (5ft. 5in.) and
stout
until emaciated by disease; he had a big head, with long black hair brushed
sideways above a big forehead; his face was square, he was clean shaven, and his
complexion never really dark, grew paler during his life in England; his ears were
small, his nose broad, and always his shining eyes were the conspicuous feature
that Ramachandra Rao observed in 1910. He walked stiffly, with head erect and toes
out-turned; if he was not talking as he walked, his arms were held clear of the
body, with hands open and palms downwards. But when he talked, whether he was
walking or standing, sitting or lying down, his slender fingers were for ever
alive, as eloquent as his countenance. }\\

The above physical description of
Ramanujan was recorded by Prof. E.H. Neville [31]. Ramanujan had only one passion
in life --– mathematics. He devoted all his time to this subject and its
development. 
Quoting Prof. Neville again [31], Ramanujan had an 
{\it instinctive perfection of
manners that made him a delightful guest or companion. Success and fame left his
natural simplicity quite untouched. To his friends he was devoted beyond measure,
and he devised curiously personal ways of showing his gratitude and expressing his
affection. The wonderful mathematician was indeed a loveable man. }\\

This is in complete accord with the views of Hardy [1] on
Ramanujan: \\

\ni {\it $\cdots$ … the picture I want to present to you is that of a man who had his
peculiarities 
like other distinguished men, but a man in whose society one could take pleasure, 
with whom one could drink tea and discuss politics or mathematics; the picture in 
short, not of a wonder from the East, or an inspired idiot, or a psychological 
fraud, but of a rational human being who happened to be a great mathematician.}\\

The integrity of Ramanujan is transparent from the following statement of Hardy 
[42]: \\

\ni {\it All of Ramanujan''s manuscripts passed through my hands, and I edited them very 
carefully for publication. The earlier ones I wrote completely. I had no share of 
any kind in the results, except of course when I was actually a collaborator, or 
when explicit acknowledgement was made. Ramanujan was almost absurdly 
scrupulous in his desire to acknowledge the slightest help.}\\

In a letter to a friend of Ramanujan, in September 1917, Hardy wrote [51]: \\

\ni {\it  He has been seriously ill but is now a good deal better. It is very difficult
to get 
him to take proper care of himself; if he would only do so we should have every 
hope that he would be quite well again before very long.} In this letter Hardy 
referred to his discovering that Ramanujan {\it  was not writing to his people nor 
apparently hearing from them. He was very reserved about it and it appeared to 
us that there must have been some quarrel.}  He expressed his anxiety regarding 
the trouble which might have arisen and wanted it to be cleared away.\\

Ramanujan was shy by temperament and contemplative by nature. He was a man 
with a great sense of humour. In the words of Neville [31]:\\

\ni {\it He had a fund of stories, and such was his enjoyment in telling them that in
his 
great days his irrepressible laughter often swallowed the climax of his
narrative.}\\ 

On learning, after his return to India, that the Government and the University of 
Madras were insisting on his going to Thanjavur, he punned on the word [52] 
Thanjavur --–  by breaking it into three part s 'anTh ',`'s'a'v'  `vnd' ''n
Tamil --
and quipped  that they wanted him to go to `'Than-savu-v''r', meaning thtownce of 
his death!  Later when he was shifted to Chetput, he punned on this word, `'ch''t'- 
`'put', and said that he was being taken to a place where everything wille 
`'very
quick''. He also did not like the name of  the building` 'Crynant' where he was 
stay in Chetput  stating that the `'Cr'y' in the word did not augur well and got 
himself shifted to another building `'Gometr'a' (which is the home where he 
breathed his last on April 26, 1920).\\

Ramanujan  was very affectionate towards his brothers and his mother, in 
particular. His wife recollected [53] that he knew astrology and made astrological 
predictions to some extent and that he knew he would not live beyond 34 years.  
Sometimes, he is supposed to have made predictions for others also. He told her 
[54], after his return from England, that he felt very happy when the Editor of  
{\it The Hindu}, Mr.Kasturiranga Iyengar, went to his room and partook the 
pongal\footnote{A South Indian delicacy prepared with rice, greengram, ghee,
pepper, jeera and cashew nuts.}  
prepared and served by him. In later years, Janakiammal told several who 
visited 
her that Ramanujan was confident his mathematics would provide her with funds, 
even after his death.\\

Some friends of Ramanujan have remembered [55] that Ramanujan could foresee 
events in visions; that being an ardent devotee of Lord Narasimha he saw drops of 
blood in dreams (which was considered as a sign of the Lord''s grace) and that  
{\it after seeing such drops, scrolls containing the most complicated mathematics 
used to unfold before him, and these he set down on paper on waking only a 
fraction of what was thus shown to him.}\\

Ramanujan''s maternal grandmother was a staunch devotee of Goddess Namagiri 
of  Namakkal. Ramanujan himself was known to his friends to be a devotee of the 
Goddess of Namakkal and he used to say that the Goddess appeared in his dreams 
and inspired him to come forth with new formulae\footnote{This was probably his 
way of explaining away his incomparable intuition and success, to those 
who could not comprehend his ability to churn out continuously new results but who
persisted in questioning him as to how he arrived at those results!}.\\
 
 Prof. K. Ananda Rao was at King's College, when Ramanujan was at Trinity 
College, and he recalled [56], in 1962, that: \\

\ni {\it In his nature he was simple, entirely free from affectation, with no trace
whatever 
of his being self-conscious of his abilities. He was quite sociable, very polite and 
considerate to others.}\\

Ramanujan never forgot that as a first born he had to shoulder the responsibility 
of taking care of his parents. He was compassionate. Accepting the University's 
offer of a scholarship, he wrote to Mr. Francis Dewsbury, the Registrar of the 
University of Madras, in a letter [57] dated 11th January 1919, from a nursing 
home in Putney:\\

\ni {\it I feel, however, that after my return to India, which I expect to happen
as
soon as 
arrangements can be made, the total amount of money to which I  shall be entitled 
will be much more than I shall require. I  should hope that, after my expenses in 
England have been paid, £ 50 a year will be paid to my parents and that the 
surplus, after my necessary expenses are met, should be used for some 
educational purpose, such in particular as the reduction of school-fees for poor 
boys and orphans and provision of books in schools. No doubt it will be possible 
to make an arrangement about this after my return. \\

\ni I feel very sorry that, as I have not been well, I have not been able to do so
much 
mathematics during the last two years as before. I hope that I shall soon be able 
to do more and will certainly do my best to deserve the help that has been given to 
me.}\\

Ramanujan concluded a letter [58] to Mr. Narayana Iyer, in November 1915, 
with the following words of gratitude:\\

\ni {\it  I am ever indebted to you and Sir Francis Spring for your zealous interest in
my case from the very beginning of acquaintance.}\\

I would like to coclude this lecture with the following assessments of
Ramanujan and his work (Bruce Berndt [60]):\\

%\leftskip = 1 cm \rightskip = 1 cm

\ni $\bullet$ {\it In notes left by B.M. Wilson, he
tells us how George Polya was captivated by Ramanujan's formulas. One day in 1951
while Polya was visiting Oxford, he borrowed from Hardy his copy of Ramanujan's 
notebooks. A couple of days later, Polya returned
them in almost a state of panic explaining that however long
he kept them, he would have to keep attempting to verify the formulae 
therein and never again would have time
to establish another original result of his own}.

\leftskip = 0 in \rightskip = 0 in

Neville began a broadcast in Hindustani, in 1941, with the declaration~: \\

%\leftskip = 1 cm \rightskip = 1 cm 

\ni $\bullet$ {\it Srinivasa Ramanujan was a mathematician so great that his name 
transcends jealousies, the one superlatively great
mathematician whom India has produced in the last thousand years}.\\ 

\leftskip = 0 cm \rightskip = 0 cm

Commenting on the quality of the theorem's in the `Lost' Notebook, Richard Askey 
says:\\

%\leftskip 1cm \rightskip 1cm

\ni $\bullet$ {\it Try to imagine the quality of Ramanujan's mind, one which drove 
him to 
work unceasingly while deathly ill, and one great enough to grow deeper while his
body became weaker. I stand in awe of his achievements; understanding is beyond 
me. We would admire any mathematician whose life's work is half of what Ramanujan
found in the last year while he was dying.}\\

\leftskip 0cm \rightskip 0cm

\ni $\bullet$ {\it Paul Erd\"os has passed on to us Hardy's personal
ratings of mathematicians: Suppose that we rate mathematicians on the
basis of pure talent on a scale from 0 to 100,
Hardy gave himself a score of 25, 
Littlewood 30, Hilbert 80, and Ramanujan 100.} (Berndt [60]).\\

\eject

\ni {\bf References}\\

\begin{enumerate} 

\item {\it Ramanujan: Twelve Lectures on subjects suggested by
his life and work}, G.H. Hardy, Chelsea, New York, 1940. % ref.1

\item {\it Collected Papers by Srinivasa Ramanujan}, 
edited by G.H. Hardy, P.V.  Seshu Aiyar and B.M. Wilson,
Chelsea, New York, 1962; first published by Cambridge Univ. Press, 1927.%ref.2 

\item {\it Ramanujan: Letters and
Reminiscences}, Memorial Number, Vol.I, ed. P.K. Srinivasan, Muthialpet High
School, Madras, 1968. % ref.3 

\item K.S. Viswanatha Sastri, in [3], P.89-93. % ref.4 

\item N. Govindarajan, in [3] P.104-–105.%ref.5
 
\item See [3] P.94, 95, 120, 121. % ref.6 

\item {\it Srinivasa Ramanujan: A Mathematical Genius}, K. Srinivasa Rao, East 
West Books (Madras) Pvt. Ltd., 1998. % ref.7 

\item {\it The Man Who Knew Infinity: A Life of the Genius Ramanujan},
Robert Kanigel, Charles Scribner's Sons, New York (1991);  Indian edition: Rupa \&
Co. (1994). % ref.8 

\item {\it Ramanujan : The Man and the Mathematician}, S.R.
Ranganathan, Asia Publishing House, 1967. % ref.9 

\item K. Chengalvarayan, in [9] p.64 (MD2). % ref.10 

\item C.R. KrishnaswamiAyyar, in [9] p. 69 (MF63). % ref.11
 
\item T. Srinivasa Raghavacharya, in [9] p. 75 (MK2). % ref.12 
 
\item R. Radhakrishna Ayyar, in [9] p. 74 (MJ91). % ref.13
 
\item N. Hari Rao, in [3] p.120-123. % ref.14 

\item V. Ramaswamy Iyer, in [3]  p.129. % ref. 15 

\item R. Krishna Rao, cousin of the mother of Prof. K. Ananda Rao. % ref. 16
 
\item R. Ramachandra Rao, in [3] p.126-127. % ref. 17

\item P.V. Seshu Aiyer in [3] p. 125. % ref. 18

\item C.L.T. Griffith to Sir Francis Spring, in [3] p.50. % ref. 19 

\item M.J.M. Hill to C.L.T. Griffith, in [3] p.53. % ref. 20

\item {\it Ramanujan: Letters and Commentary}, ed. By Bruce C. Berndt and Robert     
A. Rankin, American Mathematical Society and London Mathematical  
Society (1995); also Indian Edition with a Preface, Additions to the Indian 
Edition and Errata, by K. Srinivasa Rao, Affiliated East West Press Pvt. Ltd. 
(1997). % ref. 21

\item Ramanujan''s second letter to G.H. Hardy, in ref. 2, p. xxvii.% ref.22

\item Ref. [21], p. 17. % ref.23

\item C.P. Snow in his Forward to G.H. Hardy's{\it A Mathematician's   
        Apology}, Cambridge University Press (1976), p.30. % ref.24

\item According to C.P. Snow, Hardy was not the first eminent to be sent the 
Ramanujan manuscripts. There had been two before him, both English, both 
of the highest professional standard.  They had each returned the 
manuscripts without comment.  I don't think history relates what they said, if 
anything, when Ramanujan became famous.  As for their identity, Snow adds 
that: out of chivalry Hardy concealed this in all that he said or wrote about 
Ramanujan.  (p.33 of [24]).  However, the names are given by A. Nandy 
(in Alternative Sciences,  Allied Publishers, New Delhi, 1980) who claims 
the two to be H.F. Baker and E.W. Hobson.(also see [3] p.3). % ref.25

\item C.P. Snow in his Rectorial Address delivered before the University of St.     
        Andrews, Scotland, on 13th April 1962. % ref.26

\item Ref. [1], p.9. % ref.27

\item Ref. [3], p. 157-1158 % ref.28.

\item Ref. [3], p.55. % ref.29

\item Refer Bruce C. Berndt and Robert A. Rankin:  {\it Ramanujan: Letters and    
     Commentary}, ref. [21], for these and other letters referred to. % ref.30

\item  E.H. Neville, in ref. [3], p. 138-1141 % ref.31.

\item E.H. Neville to Dewsbury, ref. [3], p. 59-660.% ref.32

\item Littlehailes to Dewsbury, in  ref. [3], p. 61-66. % ref.33

\item Sir Francis Spring to C. B. Cotterell, in ref. [3], p. 64-665 % ref.34.

\item Letter 2 to R. Krishna Rao, in ref. [3], p. 4-7. % ref.35

\item {\it  Srinivasa Ramanujan}, J.R. Newman, in Mathematics in the modern World, 
W.H. Freeman \& Co. (1968) 73-76. % ref.36

\item  Ref. [2], p. xxx. % ref. 37

\item Ref. [8], p. 226. % ref. 38

\item Letter 4 to R. Krishna Rao, in ref. [3], p. 12-119.% ref.39

\item Letter 1 to S.M. Subramanian, in ref. [3], p. 20. % ref.40

\item P.C. Mahalanobis, in ref. [3], p. 145 – 148. Also, in{\it  Ramanujan: The Man
and the Mathematician}, S.R. Ranganathan, Asia Publishing House, 1967,  p.81. 
(MN1). % ref.41

\item G.H. Hardy to Dewsbury, in ref. [3], p. 76-777 % ref. 42

\item Ref. [21], p. 137. % ref. 43

\item Ref. [2], p. 499. % ref.44

\item Ref. [8], 233. % ref.45

\item Ref. [21], p.199. % ref.46

\item R.A. Rankin, Ramanujan as a patient, Proc. Indian Acad. Sci., Math. Sci. vol. 
93 (1984) 79-1100 % ref.47.

\item G.H. Hardy to Dewsbury, in ref. [3], p. 76-777 % ref.48.

\item Janaki Ramanujan in [9]. % ref.49

\item G.H. Hardy:  {\it A Mathematician''sApology}, (with a Foreword by C.P.
Snow), Cambridge Univ. Press (1976), first published in 1967. % ref.50

\item G.H. Hardy to Subramanian, in ref. [3], p. 68-775 % ref.51.

\item Ref. [9], p.93. (N22). % ref.52

\item Janaki Ramanujan, in ref. [3], p. 159-1161 (in Tamil), p. 17-– 172 (English 
translation). % ref.53

\item Reminiscences of Janaki Ramanujan, in ref. [9], p.89-991. (MT-– MT7 % ref.54.

\item T.K. Rajagopalan, in ref. [3], p. 167 and in ref. [9], p. 87; R. Srinivasan, in
ref. [3], p. 165 – 166; R. Radhakrishna Ayyar, in ref.[4]9, p.73 % ref.55.

\item K. Ananda Rao, in ref. [3], p. 143-144. % ref.56 

\item  Copy of Ramanujan''s letter to the Registrar, University of Madras, in ref.
[9], plate 6, between pages 104-1105. Also reproduced in ref.[ ]2, p.xix % ref.57.

\item Letter to Narayana Iyer, in ref. [3], p. 32-333 % ref.58.

\item D.A.B. Young, Ramanujan''s illness,  Current Science vol.67 (1994) p .967
- 972. % ref. 59

\item {\it Ramanujan's Notebooks}, Part I, Bruce C. Berndt, Springer-verlag (1975).
% ref. 60 

%\item Richard Askey % ref. 61.

\end{enumerate}

\end{document}